\documentclass[10pt,leqno]{article}
\usepackage{latexsym}
\usepackage[all]{xy}

\usepackage{tikz}

\usepackage{amssymb} 
\usepackage{amsmath} 
\usepackage{theorem}
\usepackage{color}

\usepackage{hyperref}

\def\K{{\mathbb{K}}}
\def\A{{\mathcal{A}}}
\def\B{{\mathcal{B}}}
\def\calS{{\mathcal{S}}}
\def\G{{\mathcal{G}}}

\DeclareMathOperator{\rank}{rank}
\DeclareMathOperator{\codim}{codim}

\DeclareMathOperator{\Der}{Der}

\numberwithin{equation}{section}

\newcommand{\owari}{\hfill$\square$}

\newtheorem{theorem}{Theorem}[section]
\newtheorem{prop}[theorem]{Proposition}
\newtheorem{cor}[theorem]{Corollary}
\newtheorem{lemma}[theorem]{Lemma}
\newtheorem{define}[theorem]{Definition}

\newtheorem{example}[theorem]{Example}

\newtheorem{conj}[theorem]{Conjecture}

\title{Derivation degree sequences of non-free arrangements}
\author{
Max Wakefield
}

\begin{document}

\maketitle

\begin{abstract}

In this note we study the logarithmic derivation module of a non-free arrangement. We prove a generalized addition theorem for all arrangements. This addition theorem allows us to find various relationships between non-free arrangements, free arrangements and restriction counts. For graphic arrangements we can use these results to find a lower bound for the maximal degree generator in terms of triangles in the associated graph. We also apply these results to the case of hypersolvable arrangements where we define hyperexponents and use them to find a lower bound for their maximal degree generator.

\end{abstract}

\section{Introduction}

The module of logarithmic derivations for hyperplane arrangements and its freeness is well studied (see \cite{Y14} for a recent survey), yet there are many unknown aspects and a compelling conjecture nearly 40 years old (Terao's conjecture see \cite{OT}). Surprisingly little has been studied about the module of logarithmic vector fields for a non-free arrangement. Fortunately, this subject is recently seeing more attention in the form of almost free or nearly free arrangements, hypersurfaces, and divisors (see \cite{DS-17-1}). Dimca and Sticlaru have studied this approach with some applications to Terao's conjecture in \cite{DS18} and  \cite{DS17-2}. The aim of this note is to take a slightly different perspective. We study arbitrary non-free arrangements and try to determine certain algebraic properties using combinatorics.

Let $\K$ be a field of characteristic zero, $V=\K^\ell$, and set $S=\calS (V^*)\cong \K[x_1,\dots ,x_\ell]$. Let $\A=\{H_1,\dots ,H_n\}$ be a finite collection of hyperplanes in $V$ with $H_i=\ker (\alpha_i)$ for all $i\in \{1,\dots ,n\}$ and put $Q(\A)=\prod\limits_{i=1}^n\alpha_i$. Denote derivations on $S$ by $$\Der (S)=\{\theta\in \mathrm{Hom}_\K (S,S)| \theta (pq)=\theta(q)p+p\theta (q)\text{ for all }p,q\in S\}.$$ Then the module of logarithmic derivations on $\A$ is defined by $$D(\A)=\{\theta\in \Der (S)| \theta (Q(\A))\in Q(\A)S\}.$$ Of course the generators for $D(\A)$ are not unique, but the if minimal system of homogeneous generators are chosen then their degrees are unique and well defined (see \cite{T-80}). This sequence of non-negative integers we will call the derivation degree sequence of the arrangement $\A$. 

If the module $D(\A)$ is free over $S$ then we say $\A$ is free and the derivation degree sequence is called the exponents of $\A$. These numbers are called the exponents because in the case of a Weyl arrangement these degree's exactly coincide with the exponents of the corresponding reflection group. More generally, for any free arrangement, these exponents are exactly the roots of the characteristic polynomial (Terao's factorization theorem \cite{Terao-gene}). For free complex arrangements the exponents have an important topological meaning. They give the Betti numbers for the complex complement via the factoring of the characteristic polynomial of the associated matroid. 

The derivation degree sequence for non-free arrangements is geometrically delicate and often is not combinatorial. For example, if you take three sets of three lines in the projective plane that intersect in a point, then the derivation degree sequence will change when the three points are colinear. On the other hand, the degree sequence gives important information like the first non-trivial Betti numbers in a free resolution of the Jacobian ideal of the arrangement. Hence they are important for understanding the singularities of the arrangement. The degree sequence also gives some information about the so called logarithmic ideal of critical sets of the arrangement (see \cite{CDFV11}). The underlying theme of this work is to use the derivation degree sequence to explore the interesting relationship between the combinatorics of the arrangement and the geometry of the arrangement.

The aim of this note is find conditions under which we can predict this derivation degree sequence combinatorially. There are three topics we study in this paper. Each of these three topics are applications of a non-free generalization of the classical ``Addition Theorem'' of Terao in \cite{T-80}. In Section \ref{sec2} we present this result that is easily generalized for adding multiple new lines (generalized from the ``free'' version in \cite{ABCHT}). In the first topic, in Section \ref{sec3}, we look at arrangements that are ``missing'' a line, meaning that there is a line not in the arrangement who has special intersection properties with the arrangement. In these cases the derivations seem to ``know'' that the line is ``missing''. This is a peculiar property and seems to have implications for Terao's conjecture. The central result in this section is an equality that relates the maximal degree generator and the minimal intersection number of a ``missing'' line. 

Next, in Section \ref{sec4}, we study the maximal degree generator in the case of graphic arrangements. The main result for graphic arrangements is a bound for this largest degree in terms of the largest number of triangles that can be made in the graph. 

Third, in Section \ref{sec5}, we apply the results of Section \ref{sec3} to hypersolvable arrangements. The hypersolvable arrangements that are not free have some interesting properties. Even though they have fairly restricted combinatorial structures the derivation degree sequence is very hard to predict without actual computation. The main result on hypersolvable arrangements is a bound on this maximal degree generator in terms of the arrangements solvable filtration. 

{\bf Acknowledgments:} We would like to thank Takuro Abe for many helpful discussions including the argument for the non-free addition theorem. The author has been partially supported by the Simons Foundation and the Naval Postgraduate School.

\section{Arrangement theory}\label{sec2}

As in the introduction we let $\K$ be a field of characteristic zero,  $V=\K^\ell$, and $S=\calS (V^*)\cong \K[x_1,\dots ,x_\ell]$ the polynomial coordinate ring. We say a hyperplane arrangement $\A=\{H_1,\dots ,H_n\}$ is central if each hyperplane, $H_i$ contains the origin and in this case the associated linear forms $\alpha_i$ are homogeneous. If $\bigcap_{H\in \A}H=\{0\}$ then we say $\A$ is essential. Given an $H\in \A$ the restriction of $\A$ to $H$ is $\A^H=\{H'\cap H| H'\in \A\}$. Since $$\Der (S)\cong \bigoplus\limits_{i=0}^\ell S\frac{\partial}{\partial x_i}$$ given a derivation $\theta\in \Der (S)$ viewed as $\theta=\sum p_i\frac{\partial}{\partial x_i}$ we say its homogeneous if all the total degrees of the $p_i$ are the same, say $d$, and we write $\deg(\theta )=d$.

The intersection lattice of $\A$ is the set $L(\A)$ of all possible intersections of hyperplanes in $\A$ with order given by reverse inclusion. We set $L(\A)_i=\{X\in L(\A)|\codim (X)=i\}$. The (single variable) M\"obius function on $L(\A)$ is given by $\mu (\mathbb{K}^\ell)=1$ and for $X>\mathbb{K}^\ell$ $$\mu(X)=-\sum\limits_{Y<X}\mu(Y).$$ Using the M\"obius function we can now define the characteristic polynomial of $\A$ by $$\chi(\A,t)=\sum\limits_{X\in L(\A)}\mu (X)t^{\dim(X)}.$$ 

Fix a minimal generating set $\{\xi_1,\dots ,\xi_k\}$ for $D(\A)$ and for $1\leq i\leq k$ let $d_i=\deg ( \xi_i)$. We call the sequence $(d_1,\dots , d_k)$ the derivation degree sequence of $\A$ and we will assume that this sequence is written in increasing order. By Saito's Criterion \cite{S80} if $k=\ell$ and the arrangement is essential then $D(\A)$ is a free $S$ module and we say that $\A$ is free. One of the hallmark Theorems in arrangement theory is Terao's factorization theorem which provides a beautiful bridge between freeness and the combinatorics on $L(\A)$.

\begin{theorem}[\cite{Terao-gene}]

If $\A$ is essential and free with derivation degree sequence $(d_1,\dots, d_\ell)$ then $$\chi(\A,t)=\prod\limits_{i=1}^\ell (t-d_i).$$

\end{theorem}

A natural next question is the so called Terao's conjecture:

\begin{conj}[\cite{OT}]

For a fixed field the freeness of $\A$ depends only on the isomorphism type of $L(\A)$.

\end{conj}

The idea of this study is to remove the restriction that $\A$ is free and try to detect combinatorial information from $D(\A)$ outside of Terao's factorization theorem.

The main tool we use in this note is a generalization of the multiple addition theorem in \cite{ABCHT} for 
the non-free case. The argument was communicated by Takuro Abe.

\begin{theorem}\label{main}
Let $\A'$ be a central arrangement in $V$ with derivation degree sequence $(d_1,\dots ,d_k)$. Let $$d=\max \{ d_i| 1\leq i\leq k\}$$ be the maximal degree generator for $D(\A')$.
 Fix
$H_1,\ldots,H_q $ hyperplanes in $V$ and let
$\A:=\A' \cup \{H_1,\ldots,H_q\}$. Assume that 
\begin{itemize}
\item[(1)] 
$H_1,\ldots,H_q$ are linearly independent,
\item[(2)]
$X:=\cap_{i=1}^q H_i \not \subset H$ for any $H \in \A'$, and 
\item[(3)]
$|\A'|-|\A^{H_i}|=d$ for $i=1,\ldots,q$.
\end{itemize}
Then there exists a set of generators $\{\theta_1,\ldots,\theta_s,\varphi_1,\ldots,\varphi_t\}$
for $D(\A')$  
such that $\deg \theta_1 \le \cdots \le \deg \theta_s<d=\deg \varphi_1=\cdots=
\deg \varphi_t$, $t \ge q$ and $$\{\theta_1,\ldots,\theta_s,\alpha_{H_1}\varphi_1,\ldots,\alpha_{H_q}\varphi_q,
\varphi_{q+1},\ldots,\varphi_t\}$$ is a set of generators for $D(\A)$.
\label{main2}
\end{theorem}

\noindent\textbf{Proof}. 
Let 
$\theta_1,\ldots,\theta_s,\varphi_1,\ldots,
\varphi_t$ be a generator for $D(\A')$ satisfying the above degree condition. At this point 
we do not know whether $t \ge q$. 
By Proposition 2.9 in \cite{T-80} by Terao, there is a polynomial $b_i\ (i=1,\ldots,q)$ of 
degree $d$ (by the condition (3)) such that 
$\theta(\alpha_{H_i}) \in Sb_i$ modulo $\alpha_{H_i}$ for $\theta \in D(\A')$. 
Let $\theta_i(\alpha_{H_j}) \equiv a_{ij} b_j$ 
modulo $\alpha_{H_j}$ for $ a_{ij} \in S\ (i=1,\ldots,s, j=1,\ldots,q)$. 
Since $\deg \theta_i<d=\deg b_j$, it holds that $a_{ij}=0\ (i=1,\ldots,s,\ j=1,\ldots,q)$, which shows that 
$\theta_i \in D(\A)\ (i=1,\ldots,s)$. 

In the same way, let $\varphi_i(\alpha_{H_j}) \equiv c_{ij} b$ modulo $\alpha_{H_j}\ (i=1,\ldots,
t,\ j=1,\ldots,q,\ c_{ij} \in S)$. Since $\deg \varphi_i=d=\deg b_j$, it holds that $c_{ij} \in \K$. 
Let $C:=(c_{ij})$ be the $(t \times q)$-matrix. Assume that 
$\mbox{rank}(C)=:r<q$. 
Then after appropriate row elementary operations, we may assume that 
$$
C=\begin{pmatrix}
E_r & *\\
O & O
\end{pmatrix}
,
$$
Then $\varphi_i \in D(\A)$ for $i=r+1,\ldots,t$. 
Now consider the tangent space $D(\A')_x$ of $D(\A')$ at $x \in X \setminus \cup_{H \in \A'} H \neq \emptyset$ 
(by the condition (2)). Hence $D(\A')_x=\K^\ell$. It is clear that 
$$
D(\A')_x=\langle 
\theta_1|_x,\ldots,\theta_s|_x,\varphi_1|_x,\ldots,\varphi_t|_x
\rangle.
$$
Since $\theta_1,\ldots,\theta_s$ and $\varphi_{r+1},\ldots,\varphi_{t}$ are tangent to 
$X$ by the above arguments, it holds that 
\begin{eqnarray*}
D(\A')_x&=&\langle 
\theta_1|_x,\ldots,\theta_s|_x,\varphi_{r+1}|_x,\ldots,\varphi_{t}|_x
\rangle\oplus 
\langle 
\varphi_{1}|_x,\ldots,\varphi_r|_x
\rangle
\\
&\subset& 
T_{X,x}\oplus
\langle 
\varphi_{1}|_x,\ldots,\varphi_r|_x
\rangle.
\end{eqnarray*}
Hence $\codim X=q=\dim_\K\langle \varphi_{1}|_x,\ldots,\varphi_r|_x\rangle$ 
by the condition (1). Thus $r \ge q$, which contradicts $r<q$. Hence $r=q$ and 
$C$ is of the form 
$$
C=\begin{pmatrix}
E_q \\
O 
\end{pmatrix}.
$$

Now define $\phi_i:=\varphi_i-\sum_{j=1}^qc_{ij}\varphi_j\ (i=q+1,\ldots,t)$ and $\phi_i:=\varphi_i \ 
(i=1,\ldots,q)$. Then 
$\phi_i(\alpha_{H_j})\equiv c_{ij}b-c_{ij}b\equiv 0$ modulo $\alpha_{H_j}$ 
for $i=q+1,\ldots,t$. Hence 
$$
\langle 
\theta_1,\ldots,\theta_s,\alpha_{H_1}\phi_1,\ldots,\alpha_{H_q}\phi_{q},\phi_{q+1}, \ldots,\phi_t
\rangle \subset
D(\A).
$$
So it suffices to show the reverse inclusion. 
Let $\theta \in D(\A)$. We prove that $
\theta \in \langle 
\theta_1,\ldots,\theta_s,\alpha_{H_1}\phi_1,\ldots,\alpha_{H_q}\phi_{q},\phi_{q+1}, \ldots,\phi_t
\rangle$. Since 
$\theta \in D(\A')=
\langle 
\theta_1,\ldots,\theta_s,\phi_1,\ldots,\phi_{t-1}, \phi_t 
\rangle$, we may express 
$$
\theta=\sum_{i=1}^s f_i \theta_i+\sum_{j=1}^t g_j \phi_j\ (f_i,g_j \in S).
$$
Hence it suffices to show that $\alpha_{H_j} \mid g_j$ for $j=1,\ldots,q$. Since 
$$
D(\A' \cup \{H_j\}) \ni \theta- \sum_{i=1}^s f_i \theta_i - \sum_{k\neq j} g_k \phi_k=g_j \phi_j,
$$
it holds that $g_j\phi_j(\alpha_{H_j}) \in S\alpha_{H_j}$. Recall that $\phi_j(\alpha_{H_j})=
\varphi_j(\alpha_{H_j}) \equiv b_j$ modulo $\alpha_{H_j}$ and $b_j \not \in S\alpha_{H_j}$. Hence 
$g_jb_j \in S\alpha_{H_j}$ shows that $g_j \in S \alpha_{H_j}$. As a consequence, 
$$
\theta=\sum_{i=1}^s f_i \theta_i+\sum_{j=r+1}^{t} g_j \phi_j+
\sum_{j=1}^r (g_j/\alpha_{H_j})\alpha_{H_j} \phi_j
$$
with $g_j/\alpha_{H_j} \in S$, which completes the proof. \owari
\medskip

\begin{example}
Theorem \ref{main2} does not imply that, if the condition in Theorem \ref{main2} holds and 
$\A'$ is not free, then $\A$ is not free. For example, examine the arrangement 
$$
\A':=\{xy(x-z)(y-z)=0\}.
$$
This is not free and has a generators 
$$
\theta_E,\ x(x-z)\partial_x,\ 
y(y-z)\partial_y,\ 
(x-z)(y-z)\partial_z.
$$
Hence the derivation degree sequence is $\{1,2,2,2\}$. Hence in the setup of Theorem \ref{main2}, 
$d=2$. Let $H:=\{x-y=0\}$ and $\A:=\A' \cup \{H\}$. Then it is well-known that $\A$ is free with 
$\exp(A)=(1,2,2)$. However, the construction in Theorem \ref{main2} presents a generating set of $D(\A)$ with degrees $\{1,2,2,3\}$. In this case, the generating set from Theorem \ref{main} is  
$$
\{\theta_E,\ x(x-z)\partial_x+y(y-z)\partial_y,\ 
(x-z)(y-z)\partial_z, (x-y)y(y-z)\partial_y\}. 
$$
However, it is clear that the last derivation is not necessary to generate $D(\A)$. Namely, 
$$
\{\theta_E,\ x(x-z)\partial_x+y(y-z)\partial_y,\ 
(x-z)(y-z)\partial_z\}$$
form a basis for $D(\A)$. Hence in particular, Theorem \ref{main2} does not 
necessarily preserve the minimality of the generators.
\end{example}

\section{The Largest Degree Generator and Minimal Restrictions}\label{sec3}

In this section we examine relations between degrees of generators of $D(\A)$ and 
the number of hyperplanes in the restricted arrangement. We start with a simple fact about how restriction numbers give a lower bound on the maximal degree generator.

\begin{prop}\label{genbound}
Let $\A$ be a central arrangement, $H \in \A$ and 
$\A':=\A \setminus \{H\}$. Assume that $|\A^H| \le |\A'|-d$. Then 
every set of generators of $D(\A')$ has to contain an element $\theta \in D(\A')$ with 
$\deg \theta \ge d$. 
\label{degree}
\end{prop}

\noindent
\textbf{Proof}. 
Let $\{\theta_1,\ldots,\theta_s\}$ be a set of generators of $D(\A')$. 
Assume that $\deg \theta_i <d$ for $i=1,\ldots,s$. Then again by Proposition 2.9 in \cite{T-80} there is a polynomial $b$
such that $\theta_i(\alpha_H)\equiv c_i b$ modulo $\alpha_H$ for some $c_i$. Also, $\deg(b)=|\A'|-|\A^H|\geq d$ and $\deg(\theta_i(\alpha_H)<d$ implies
$c_1=\cdots=c_s=0$. Hence $D(\A') \subset D(\A) \subset D(\A')$, which 
is a contradiction. \owari
\medskip

Proposition \ref{degree} gives a non-freeness criterion as follows.

\begin{cor}
Let $\A$ be a central arrangement, $H \in \A$ and 
$\A':=\A \setminus \{H\}$. Assume that 
$\chi(\A',t)=\prod_{i=1}^\ell (t-d_i)$ with 
$d_1 \le \cdots \le d_\ell$. If $|\A^H| \le \sum_{i=1}^{\ell-1} d_i-1$, then 
$\A'$ is not free.
\label{nonfree}
\end{cor}

\noindent
\textbf{Proof}. 
Note that $|\A^H| \le \sum_{i=1}^{\ell-1} d_i-1=|\A'|-(d_\ell+1)$. 
Assume that $\theta_1,\ldots,\theta_\ell$ form a free basis for $D(\A')$ with 
$\deg \theta_i=d_i$. Then Corollary \ref{degree} shows that 
$\max_{i=1}^\ell \{d_i\}=d_\ell \ge d_\ell+1$, which is a contradiction. \owari
\medskip

Corollary \ref{nonfree} gives a nice combinatorial method to show non-freeness in certain examples where the characteristic polynomial factors. For example, to our knowledge the following example presents a new combinatorial method, which does not reply on the addition-deletion theorem, to show the arrangement of Example 4.139 in \cite{OT} is not free.

\begin{example}\label{star+}
Let $\A':=\{x(x\pm y)(x\pm 2y)(y-z)z=0\}$. Then 
$H:=\{y=0\}$ intersects with $\A'$ at $2$-points. 
Since $\chi(\A',t)=(t-1)(t-3)^2$ and 
$2=|(\A' \cup \{H\})^H|\le 1+3-1=3$, Corollary \ref{nonfree} shows that 
$\A'$ is not free. 
\end{example}

Another application is a variant of the addition-type theorem from 
non-free to free arrangements. The arrangements that satisfy the hypothesis of the following Proposition seem to be close to nearly free arrangements.

\begin{prop}\label{4gens}
Let $\A'$ be a arrangement in $\mathbb{K}^3$ such that 
$D(\A')$ has a minimal generators $\theta_E,\theta_1,\theta_2,\theta_3$ such that 
$\deg \theta_i=d_i$ for $i=1,2,3$ and that $1 \le d_1 \le d_2 \le d_3$. If there is a plane 
$H \not \in \A'$ such that $|\A'|-|\A^H|=d_3$ for $\A:=\A' \cup \{H\}$ and 
$1+d_1+d_2=|\A|$, then 
$\A$ is free with $\exp(\A)=(1,d_1,d_2)$.
\label{free}
\end{prop}

\noindent
\textbf{Proof}. 
By Theorem \ref{main2}, we may assume that 
$\theta_1,\theta_2 \in D(\A)$ and $\theta_3 \not \in D(\A)$. If 
$\theta_E,\theta_1,\theta_2$ are linearly independent over $S$, then 
Saito's criterion completes the proof. So assume that they are not 
$S$-independent. Since $D(\A)$ has a resolution 
$$
0 \rightarrow S \rightarrow S^4 \rightarrow D(\A) \rightarrow 0,
$$
it holds that $D(\A)=M \oplus S \theta_3$ for 
$M:=\langle \theta_E,\theta_1,\theta_2 \rangle_S$. Let 
$M_0:=\{\theta \in M \mid \theta(Q)=0\}$ for $Q:=Q(\A)$. 
Then we show that $M =M_0 \oplus S\theta_E$. 

Define a map $f:M \rightarrow M_0$ by 
$f(\theta):=\theta-\theta(Q)\theta_E/(\deg Q)Q$. Since 
$M \subset D(\A)$, $f$ is defined and surjective. 
Also, a canonical inclusion $M_0 \subset M$ gives a section of 
$f$. Since 
$\mbox{ker}(f)=S \theta_E$, it holds that 
$M=M_0 \oplus S \theta_E$ as required. 

Hence the linear relation among $\theta_E,\theta_1,\theta_2,\theta_3$ is in $M_0$. 
Hence we may replace $\theta_1$ and $\theta_2$ by $\varphi_1$ and $\varphi_2$ respectively 
such that $\varphi_i \in D(\A)$, $\theta_E,\varphi_1,\varphi_2,\theta_3$ form a generator 
for $D(\A')$ and $\langle \varphi_1,\varphi_2\rangle_S=M_0$. Hence there is a 
polynomial $g,h \in S$ such that $g\varphi_1=h \varphi_2$. We may assume that $g$ and $h$ are 
coprime. Hence $\varphi_1/h=\varphi_2/g=:\phi$ is a regular derivation. 
Now we want to show that $\phi \in D(\A)$. Assume $\phi \notin D(\A)$. Then there is a plane $L \in \A$ such that $\phi(\alpha_L) \not \in S \alpha_L$. 
Since $h\phi=\varphi_1$ and $g\phi=\varphi_2$ are both in $D(\A)$, we know that 
both $h\phi(\alpha_L)$ and $g \phi(\alpha_L)$ are divisible by $\alpha_L$. Since 
$\phi(\alpha_L)$ is not divisible by $\alpha_L$ and $g$ and $h$ are coprime, this is a contradiction. 
Hence $\phi \in D(\A)$. 

Then $D(\A') \supset \langle \theta_E,\phi,\theta_3\rangle \supset 
\langle \theta_E,\varphi_1,\varphi_2,\theta_3\rangle =D(\A')$ contradicts the minimality of 
$\theta_E,\theta_1,\theta_2,\theta_3$ as a generators, which completes the proof. \owari

\

We also give an another application for arrangements in $\mathbb{K}^3$.

\begin{cor}
Let $\A'$ be a central arrangement in $\mathbb{K}^3$, $H \not \in \A'$ a hyperplane and 
$\A:=\A' \cup \{H\}$ with $\chi(\A,t)=(t-1)(t-a)(t-b)\ (a \le b)$. Assume that 
the derivation degree sequence of $\A$ is
$(d_1,\dots ,d_k)$ where $k \ge 5$. If $|\A'|-|\A^H| =d_k$ and 
$b \le d_k$, then $\A$ is not free.
\label{3nonfree}
\end{cor}

\noindent
\textbf{Proof}. 
By Theorem \ref{main2}, we may assume that 
there is a set of minimal generators $\{\theta_1,
\ldots,\theta_{s},\varphi_1,\dots ,\varphi_t\} $ for $D(\A')$ such that $\deg (\varphi_i)=d_k$ is the maximal degree and $\{\theta_1,
\ldots,\theta_{s},\alpha_H\varphi_1,\varphi_2,\dots ,\varphi_t\}$ generate $D(\A)$. Suppose that $\A$ is free. Since $a\leq b$ and $d_k\geq b$ by Terao's factorization theorem we know that $\alpha_H\varphi_1$ will not be a minimal generator. Hence we may choose a set of 3 minimal generators for $D(\A)$ from  $\{\theta_1,
\ldots,\theta_{s},\varphi_2,\dots ,\varphi_t\}$ and $s+t-1=k-1\geq 4$. This contradicts the minimality of the set of generators for $D(\A')$. \owari

\

The next proposition provides bounds on how much the number of generators of certain degrees can increase when adding a hyperplane. For a central arrangement $\B$ with ordered degree sequence $(d_1,\ldots ,d_k)$, let $n_i(\B)=\max \{j| d_j\leq i\}$.

\begin{prop}
Let $\A'$ be a central arrangement and $\A:=\A' \cup\{H\}$. Suppose that the maximal degree of a minimal set of generators for 
$D(\A')$ is $d$.
\begin{enumerate}
\item If $|\A'|-|\A^H|=d$ then $n_d(\A') -1 \le n_d(\A)$.

\item If $|\A'|-|\A^H| =e+1 \le d$ then $n_e(\A')  \le n_e(\A)$.
\end{enumerate}
\label{numgen}
\end{prop}

\noindent
\textbf{Proof}.
1.\,\,
Let $\theta_1,\ldots,\theta_s,\varphi_1,\ldots,\varphi_t$ be a minimal set of generators for $D(\A')$ such that $\deg \theta_1=:d_1 \le 
\cdots \le d_s:=\deg \theta_s <d:=\deg \varphi_1=\cdots=\deg \varphi_t$. Any other minimal set of generators for $D(\A')$ consists of $(s+t)$-derivations. 
By applying Theorem \ref{main2}, we may assume that $\theta_1,\ldots,\theta_s,\varphi_1,\ldots,\varphi_{t-1}, 
\alpha_H \varphi_t \in D(\A)$. Assume that $n_d(\A) \le n_d(\A')-2=s+t-2$. 
Let $\phi_1,\ldots,\phi_u$ be a part of the minimal-number generators for $D(\A)$ the degrees of which 
are at most $d$ and $u \le s+t-2$. Since $\langle \phi_1,\ldots,\phi_u \rangle \supset 
\langle \theta_1,\ldots,\theta_s,\varphi_1,\ldots,\varphi_{t-1}\rangle$, we may choose 
$\phi_1,\ldots,\phi_u,\varphi_t$ as a generator for $D(\A')$. Since 
$u+1 <s+t$, this contradicts the minimality. 

2. \,\,
Let $\theta_1,\ldots,\theta_s,\varphi_1,\ldots,\varphi_t$ be a minimal set of generators for $D(\A')$ such that $\deg \theta_1=:d_1 \le 
\cdots \le d_s=\deg \theta_s <e+1=\deg \varphi_1\le \cdots \le \deg \varphi_t=d$. 
Assume that $n_e(\A')>n_e(\A)$. By the arguments in the proof of Theorem \ref{main2}, 
we know that $\theta_1,\ldots,\theta_s \in D(\A)_{\le e}$. Let 
$\phi_1,\ldots,\phi_u$ be the part of a minimal set of generators for $D(\A)$ such that 
their degrees are at most $e$. Then the assumption says that $u <s$. Then 
$\langle \phi_1,\ldots,\phi_u\rangle \supset \langle \theta_1,\cdots,\theta_s\rangle$. Thus 
$\phi_1,\ldots,\phi_u,\varphi_1,\ldots,\varphi_t$ form a generator for $D(\A')$, which contradicts the minimality 
of the original generating set.\owari
\medskip

The next proposition begins our discussion of bounding restriction numbers with degrees of generators.

\begin{prop}\label{proprestrict1}

Let $\A'$ be a central arrangement in $\K^\ell$ with ordered derivation degree sequence $(d_1,\dots ,d_k)$. Then there does not exist any hyperplane $H\subset V$ with $H\not\in \A'$ such that $|\A'|-d_k>|\A^H|$.

\end{prop}

\textbf{Proof}. As in Theorem \ref{main2} suppose $\{\theta_1,\ldots,\theta_s,\varphi_1,\ldots,\varphi_t\}$ is a minimal set of generators of $D(\A')$ where $d_k=\deg \varphi_1=\cdots=\deg \varphi_t$. Suppose that there does exist an $H\not\in \A$ such that $|\A'|-d_k>|\A^H|$. Then $|\A'|-|\A^H|>d_k$. Again from Proposition 2.9 in \cite{T-80} there exists a polynomial $b$ such that $\deg b>d_k$ and $\phi_i(\alpha_H)\in (\alpha_H,b) $ for all $i$. Since $\deg \phi_i(\alpha_H)<d_k$ we have that $\phi_i\in D(\A)$ which is the same reasoning as that given in the proof of Theorem \ref{main2} that $\theta_i\in D(\A)$. Hence $D(\A')\subseteq D(\A)$ but by definition $D(\A)\subseteq D(\A')$. This implies that $H\in \A'$ which is a contradiction. \owari

\

The next corollary follows directly.

\begin{cor}\label{countpoints}

Let $\A$ be a central arrangement in $\mathbb{K}^3$ with ordered derivation degree sequence $(d_1,\ldots ,d_k)$. Let $H$ be any plane in $\mathbb{K}^3$ that is not in $\A$. Then $|H\cap\bigcup\A|\geq |\A|-d_k$.

\end{cor}

In order to make these results a little more clear we propose the following notation. We call $$t_\A=\min \left\{ \left|(\A\bigcup H)^H\right| :H \text{ is any hyperplane in } \mathbb{K}^\ell \text{ not in }\A \right\} $$ the \emph{minimal restriction number} of $\A$ and let $d_\A=d_k$ where $(d_1\ldots,d_k)$ is an ordered derivation degree sequence for $\A$. Then we can rephrase Proposition \ref{proprestrict1} and Corollary \ref{countpoints} as \begin{equation}\label{unequal} t_\A\geq |\A | -d_\A .\end{equation}
We are interested in which arrangements have $t_\A=|\A |-d_\A$ because arrangements with this property seem to have strong connections between the derivation module and their combinatorics. The following example is just a generalization of Example \ref{star+}, but it does give an infinite family where $t_\A=|\A |-d_\A$.

\begin{example}

Let $\B_n$ be defined by the product $$xz(y-z)\prod\limits_{k=1}^n(x+ky) .$$ Note that Example \ref{star+} is a projective transformation of $\A_4$. Then the module of derivations $D(\B_n)$ is generated by $$\{\theta_E, z(y-z)\partial_z,x\prod\limits_{k=1}^n(x+ky)\partial_x,(y-z)\prod\limits_{k=1}^n(x+ky)\partial_y\} .$$ Here $t_{\B_n}=2$, $|\B_n|=n+3$, and $d_{\B_n}=n+1$.

\end{example}

Of course, the inequality \ref{unequal} is in some sense rarely an equality. The next example is an infinite family of arrangements where $t_\A> |\A | -d_\A$.
\begin{example}

Let $\G_{2,n}$ be a generic arrangement of $n>3$ hyperplanes in $\mathbb{K}^3$ where the hyperplanes are in general position. In \cite{Yuz91} Yuzvinsky proved that $D(\G_{2,n})$ is minimally generated by $\theta_E$ and other derivations of degree $n-2$. The smallest that $|\G_{2 ,n}^H|$ could be is $n-2=t_{\G_{2,n}}$. But $|\G_{2,n}|-d=n-(n-2)=2$. So equality only happens when $n=4$.

\end{example}

The next example shows that equality of \ref{unequal} does not hold in general even for free arrangements.

\begin{example}

Let $\A_3$ be the braid arrangement defined by $Q=(x-y)(x+y)(x-z)(x+z)(y-z)(y+z)$ in $\mathbb{C}^3$. Then $t_{\A_3}=4$  (actually this number will be further justified in the next section) and  $|\A_3|-d_\A=6-3=3$.

\end{example}

Clearly the minimal restriction number $t_\A$ is not an invariant of the intersection lattice. However, some arrangements have the property that $t_\A$ is determined combinatorially. The following is such a class.

\begin{prop}\label{2-points}

Suppose that $\A$ is an arrangement in $\mathbb{K}^3$ such that there are two intersections $p_1,p_2\in L(\A)_2$ such that $\mu(p_1)+\mu(p_2)=|\A|-2$ and there is not a line in $\A$ that contains both $p_1$ and $p_2$. Then $t_\A=2$.

\end{prop}

With these types of arrangements using \ref{unequal} we get a large bound on the largest degree minimal generator.

\begin{cor}

For arrangements that satisfy the conditions of Proposition \ref{2-points} $$d_\A\geq |\A|-2.$$

\end{cor}

\begin{example}

Let $\A$ be defined by the polynomial $Q=yz(x-z)(x+z)(x-y)(x+y)$. Then $\A$ satisfies the conditions of Proposition \ref{2-points} and $t_\A=2$. The degrees of a minimal generating set for $D(\A)$ is $(1,3,3,4)$ thanks to Macaualy2. Hence $d_\A=4=|\A|-t_\A$.

\end{example}



\section{Graphic Arrangements}\label{sec4}

Let $G=(V,E)$ be an undirected simple graph. Then the graphic arrangement corresponding to $G$ is $$\A_G=\{\{x_i-x_j=0\}| (i,j)\in E\}$$ which is a sub arrangement of the braid arrangement $\A_{|V|}$. In this section we will compute the minimal restriction $t_{\A_G}$ for graphic arrangements in terms of the graph. We will need a few trivial lemmas.

\begin{lemma}\label{int1}

Suppose that $X\in L(\A_G)_2$ comes from the intersection of two edges in $G$ which share a point. Without loss of generality we may assume that the edges are $(1,2)$ and $(2,3)$ so that $X=\{(x,x,x,z_4,\dots ,z_\ell)|x,z_i\in \mathbb{K}\}$. If any hyperplane $H$ with linear form $\alpha_H=\sum_{i=1}^\ell c_ix_i$ contains $X$ then $c_1+c_2+c_3=0$ and $c_i=0$ for $4\leq i\leq \ell$.

\end{lemma}



\begin{lemma}\label{int2}

Suppose that $X\in L(\A_G)_2$ comes from the intersection of two edges in $G$ which do not share a point. Without loss of generality we may assume that the edges are $(1,2)$ and $(3,4)$ so that $X=\{(x,x,y,y,z_5,\dots ,z_\ell)|x,y,z_i\in \mathbb{K}\}$. If any hyperplane $H$ with linear form $\alpha_H=\sum_{i=1}^\ell c_ix_i$ contains $X$ then $c_1+c_2=0$, $c_3+c_4=0$, and $c_i=0$ for $5\leq i\leq \ell$.

\end{lemma}

The next lemma is the crux of the computation of $t_{\A_G}$. 

\begin{lemma}\label{cont}

Let $X,Y\in L(\A_G)_2$ with $X\neq Y$. There exists a hyperplane $H\not\in \A_G$ that contains both $X$ and $Y$ if and only if the edges of $X$ and $Y$ make a 4-cycle. Moreover, this hyperplane is unique. 

\end{lemma}

\textbf{Proof}. Suppose that a hyperplane $H$ contains two different flats $X,Y\in L(\A_G)_2$. Note that if $X$ and $Y$ come from intersections of edges that share an edge then the only hyperplane that would contain both $X$ and $Y$ would be the hyperplane given by that shared edge. The rest of the proof is broken up into cases.

First look at the case where both $X$ and $Y$ come from two edges that share a common point. Suppose that $X=\{(z_0,\dots ,z_\ell)| z_i=z_j=z_k \}$ and $Y=\{(z_0,\dots ,z_\ell)| z_{i'}=z_{j'}=z_{k'} \}$. If a hyperplane $H$ with corresponding linear form $\alpha_H=\sum_{i=0}^\ell c_ix_i$ to contains $X$ and $Y$ then exactly two of the vertices must be the same other wise Lemmas \ref{int1} would not be satisfied. Without loss of generality we may assume the two shared vertices are $i$ and $j$, and so by Lemmas \ref{int1} again we get that $\alpha_H \doteq x_i-x_{j}$. Since $H$ is not in $\A$ and $X$ and $Y$ do not share an edge we know that the vertices $\{i,j,k,k'\}$ must form a 4-cycle. 

For the second case, suppose that $X$ comes from the intersection of two edges that share a point and $Y$ comes from the intersection of two edges that do not share a point. Then we can say that $X=\{(z_0,\dots ,z_\ell)| z_i=z_j=z_k \}$ and $Y=\{(z_0,\dots ,z_\ell)| z_{a}=z_{b} \text{ and }z_c=z_{d} \}$. At least two of the vertices from $\{i,j,k\}$ must match $\{a,b,c,d\}$ other wise by Lemmas \ref{int1} and \ref{int2} the coefficients of the hyperplanes linear form will be zero. WLOG assume that $j$ is the middle vertex in $\{i,jk\}$. There are 4 cases to check: (1) $i=a$, $k=b$, $c,d\neq j$ (2) $i=a$, $k=b$, $c=j$, and $d\neq i,j,k$ (3) $i=a$, $c=k$, $b\neq j,k$, and $d\neq i,j$ (4) $i=a$, $j=c$, $b\neq i,j,k$, $d\neq i,j,k$. In each cases (1) and (2) the only hyperplane that contains $X$ and $Y$ is $x_i-x_k=x_a-x_b$ which is in $\A$. In cases (3) and (4) the restrictions imply the coefficients of the hyperplane that contains both $X$ and $Y$ must all be zero.

The third case is when both $X$ and $Y$ come from intersections of edges that do not share a common vertex. Then we can say that $X=\{(z_0,\dots ,z_\ell)| z_a=z_b \text{ and } z_c=z_d \}$ and $Y=\{(z_0,\dots ,z_\ell)| z_{a'}=z_{b'} \text{ and }z_{c'}=z_{d'} \}$. From Lemma \ref{int2} if a hyperplane $H$ contains $X$ then $\alpha_H=s(x_a-x_b)+t(x_c-x_d)$ where $s$ and $t$ are non-zero constants. The only case where $H$ could vanish on $Y$ without sharing an edge is when $\{a,b,c,d\}=\{a',b',c',d'\}$  and the union of the two sets of edges make a 4-cycle. Now there are two cases, (1) $a=a'$, $b=c'$, $c=b'$ and $d=d'$ and (2) $a=a'$, $d=b'$, $b=b'$, and $c=c'$. In case (1) the only hyperplane that would vanish on both $X$ and $Y$ is $\alpha_H\doteq x_a-x_b+x_c-x_d$. In case (2) the only hyperplane that would vanish on both $X$ and $Y$ is $\alpha_H\doteq x_a-x_b-x_c+x_d$. \owari

Let $Tri(G)=$maximal number of new triangles that can be made from adding an edge to $G$.

\begin{lemma}\label{codim2max}

Let $r_{\A_G}=$ maximal number of codimension 2 flats that any hyperplane $H$ not in $\A_G$ can contain of $L(\A_G)$. Then$$r_{\A_G}=\left\{ \begin{array}{lcl}
0&\hspace{.5cm}& \text{ if } |E|=1\\
1 &&\nexists\text{ 4-cycle and } |E|\geq 2\\
2 &&  Tri(G)=0\text{ and } \exists \text{ 4-cycle}\\
Tri(G) && \text{otherwise}.

\end{array}\right. $$

\end{lemma}

{\bf Proof}. The case $|E|=1$ is trivial. If there are no 4-cycles then by Lemma \ref{cont} then no hyperplane can contain two non-equal codimension 2 flats.  Now suppose that $Tri=0$ and there is a 4-cycle. Let $X$ be two opposing sides of the 4-cycle and $Y$ be the opposite two opposing sides. From the third case in the proof of Lemma \ref{cont} we know there is a hyperplane which contains both $X$ and $Y$. Hence in this case $r_{\A_G}\geq 2$. Now suppose there is a hyperplane $H$ that contains 3 codimension 2 intersections $X,Y,Z\in L(\A_G)$. Since no new triangles can be formed we know that we can not have the second case in the proof of Lemma \ref{cont} because in that case the uniquely defined hyperplane that contains the subspaces is of the type $x_i-x_j$. Hence the edges of $X$, $Y$, and $Z$ must pairwise form 4-cycles. This can only happen if the subgraph is $K_4$ the complete graph on 4 vertices. This is covered in the third case of the proof of Lemma \ref{cont} and it is noted that for the two different types of 4-cycles two distinct hyperplanes uniquely contain the two respective pairs of subspaces. Hence there can not be 

Now suppose that $Tri>0$. This is the second case in the proof of Lemma \ref{cont}. For each new triangle that can be made the hyperplane which completes the triangle will be containing the corresponding subspace for the other two edges. Hence $r_{\A_G}=Tri(G)$. \owari

\

The next theorem follows directly from Lemma \ref{codim2max} since the number $t_{\A}$ for any arrangement is the size of the arrangement minus the maximal number of codimension 2 intersections that a hyperplane not in $\A$ can contain i.e. $t_{\A}=|\A|-r_\A$.

\begin{theorem}\label{triangles}

$$t_{\A_G}=\left\{ \begin{array}{lcl}
|E|&\hspace{.5cm}& \text{ if } |E|=1\\
|E|-1 &&\nexists \text{ 4-cycle and } |E|\geq 2\\
|E|-2 && Tri(G)=0\text{ and } \exists \text{ 4-cycle}\\
|E|-Tri(G) && \text{otherwise}.

\end{array}\right. $$
\end{theorem}

\

Using Theorem \ref{triangles} and (\ref{unequal}) we get the following corollary which presents a nice lower bound on the largest degree generator.

\begin{cor}

If $G$ is a graph where new triangles can be made then $$d_{\A_G}\geq Tri(G).$$

\end{cor}

Now we present an example where the bound on $d_{\A_G}$ is tight.

\begin{example}

Let $G=(\{1,\dots ,8\},E)$ where $$E=\{(1,3),(1,4),(1,5),(1,6),(1,7),(1,8),(2,3),(2,4),(2,5),(2,6),(2,7),(2,8)\} .$$ Then $Tri(G)=6$ and a Macaulay 2 calculation shows that a minimal generating set for $D(\A_G)$ has degrees $(1,2,2,2,2,2,2,6)$. So $d_{\A_G}=Tri(G)$. 

\end{example}

The next example shows that there are examples where $Tri(G)>0$ and $d_{\A_G}>Tri(G)$.

\begin{example}

Let $G=(\{1,2,3,4,5\},E)$ where $$E=\{(1,2),(1,3),(1,4),(2,3),(2,4),(3,4),(3,5),(4,5)\} .$$ $Tri(G)=2$ but $d_{\A_G}=3$. 

\end{example}

In this last example, $K_4$ is a subgraph and the more general inequality $d_{\A_G}\geq r_{\A_G}$ is not tight with $G=K_4$. {\bf Question:} Is there a nice formula for $d_{\A_G}$ in terms of $Tri(G)$ and the maximal complete subgraph in $G$?

\section{Hypersolvable Arrangements}\label{sec5}

The goal for this section is to acquire some information about the degrees of generators for hypersolvable arrangements using their combinatorial structure. As in the case of graphic arrangements we focus on the maximal degree generator. Hypersolvable arrangements were first defined by Jambu and Papadima in \cite{JP98} and studied from a topological perspective. They showed that the $K(\pi,1)$ property was combinatorial within the class of hypersolvable arrangements, they gave presentations for their fundamental groups, and they showed that the quadratic Orlik-Solomon algebra of a hypersolvable arrangement is Koszul. We now review the definition. Let $\A $ and $\B$ be arrangements such that $\rank (\A)>2$ and $\B\subseteq \A$. Set $\bar{\B}=\A\backslash \B$. For elements $H_{i_1},\dots ,H_{i_k}\in \A$ we define $\rank (H_{i_1},\dots ,H_{i_k})=\codim (\bigcap H_{i_j})$.

\begin{define}

We say that $\B$ is \emph{solvable} in $\A$ if 

\begin{enumerate}

\item\label{it1} For all $\alpha,\beta \in \B$ with $\alpha\neq\beta$ and for all $a\in \bar{\B}$ the $\rank (\alpha,\beta,a)=3$.

\item\label{it2} For all $a,b\in \bar{\B}$ there exists $\alpha\in \B$ such that $\rank (a,b,\alpha)=2$. By \ref{it1}. we have that $\alpha $ is unique, hence we will call it $f(a,b)$.

\item For all $a,b,c\in \A$ $\rank (f(a,b),f(a,c),f(b,c))=1\text{ or }2$.

\end{enumerate}

\end{define}

Now we define an arrangement to be hypersolvable if it has a chain of solvable subarrangements.

\begin{define}

An arrangement $\A$ is \emph{hypersolvable} if there exists subarrangements $\emptyset=\A_0\subset \A_1\subset \A_2 \subset \cdots \subset \A_k=\A$ such that for all $i=1,\dots , k-1$ the subarrangement $\A_i$ is solvable in $\A_{i+1}$. We will call this filtration a \emph{solvable filtration}.

\end{define}

Using \cite[Proposition 1.10]{JP98} we can define supersolvable arrangements as follows.

\begin{define}

An arrangement $\A$ is \emph{supersolvable} if it is hypersolvable with $k=\rank(\A)$.

\end{define}

It is well known that supersolvable arrangements are free. Here we recall how to find the degrees of the generators of the module of derivations by rephrasing \cite[Theorem 4.58]{OT}.

\begin{theorem}\label{superthm}

Let $\A$ be a supersolvable arrangement with a solvable filtration $\emptyset=\A_0\subset \A_1\subset \A_2 \subset \cdots \subset \A_{\rank{\A}}$. Set $b_i=|\A_i\backslash \A_{i-1}|$ for $i=1,\dots ,\rank{\A}$. Then the exponents of $\A$ are $$\mathrm{exp}(\A)=(b_1,\dots ,b_{\rank{\A}}).$$

\end{theorem}

Ideally it would be interesting if there was some kind of generalization of Theorem \ref{superthm} to hypersolvable arrangements. The most natural generalization (the generator degrees are the $b_i$ from above) does not hold (see example \ref{super1ex}) but at least we can get a kind of bound for the top degree generator. To do this we make the following definition.

\begin{define}

Suppose that $\A$ is a hypersolvable arrangement with solvable filtration $\emptyset = \A_0\subset \A_1\subset \cdots \subset \A_k=\A$. The \emph{hyperexponents} of $\A$ is the multiset of positive integers $\mathrm{hypexp}(\A, \{A_i\})=\{b_1,\dots ,b_k\}$ where $b_i=|\A_i\backslash \A_{i-1}|$.

\end{define}

Recalling a result from \cite{JP98} we have that $\mathrm{hypexp}(\A,\{A_i\} )$ does not depend on the filtration. To do this we summarize some of Section 3 of \cite{JP98}. Let $R$ be a commutative ring and $\Lambda^*_R(\omega_H |H\in \A)$ be the graded commutative exterior algebra with coefficients in $R$. The Orlik-Solomon ideal of an arrangement $\A=\{H_1,\dots ,H_n\}$ is $$I(\A)=\Bigg\langle  \sum\limits_{j=1}^{j-1}(-1)^s\omega_{H_{i_1}}\cdots \hat{\omega}_{H_{i_j}}\cdots \omega_{H_{i_s}}\Bigg | \rank \{H_{i_1},\dots ,H_{i_s}\}<s \Bigg\rangle$$ where the hat symbol means to remove that element from the product. Probably the most celebrated result in the theory of hyperplane arrangements is that if $\A$ is a collection of hyperplanes in a complex vector space and $M(\A)$ the complement of the union of the hyperplanes then the \emph{Orlik-Solomon algebra} $$OS(\A,R)=\Lambda^*_R(\omega_H |H\in \A)/I(\A)$$ is isomorphic to the cohomology algebra $OS(\A,R)\cong H^*(M(\A),R)$. For our context we focus on a similar algebra which records only the rank 2 dependences. The quadratic Orlik-Solomon ideal is $$QI(\A)=\big \langle \omega_{H_i}\omega_{H_j}-\omega_{H_i}\omega_{H_k}+\omega_{H_j}\omega_{H_k} \big | \rank \{H_i,H_j,H_k\}=2\big \rangle $$ and the \emph{quadratic Orlik-Solomon algebra} is $$QOS(\A,R)= \Lambda^*_R(\omega_H |H\in \A)/QI(\A).$$ Now following \cite[Definition 3.1]{JP98} the \emph{quadratic Poincar\'e polynomial} of $\A$ is $$QP(\A,t)=\sum\limits_{i\geq 0}\dim_R (QOS(\A,R)^i) t^i$$ which does not depend on $R$. Now we can state \cite[Proposition 3.2]{JP98}.

\begin{prop}\label{quadprod}

Suppose that $\A$ is a hypersolvable arrangement with solvable filtration $\emptyset = \A_0\subset \A_1\subset \cdots \subset \A_k=\A$. Then $$QP(\A,t)=\prod\limits_{i=1}^k (1+b_it)$$ where as above $b_i=|\A_i\backslash \A_{i-1}|$.

\end{prop}

Because of Proposition \ref{quadprod} the hyperexponents do not depend on the filtration and we will denote them by $\mathrm{hypexp}(\A)$. Now we turn our focus to a few simple lemmas in order to state and prove the main result of this section.

\begin{lemma}\label{restrictsolv}

If a subarrangement $\B$ in $\A$ is solvable then for any $H\in \bar{B}$,  $|\A^H|=|\B|$.

\end{lemma}

\noindent\textbf{Proof}. By \ref{it1}. we have that $|\A^H|\geq |\B|$ and by \ref{it2}. we have that $|\A^H|\leq |\B|$. \owari

\begin{lemma}\label{maxup}

Let $\B\subseteq \A$ be arrangements with maximal degree generators $d_\B$ and $d_\A$ respectively. Then $d_\B\leq d_\A$.

\end{lemma}

\noindent\textbf{Proof}. Suppose that $d_\B> d_\A$ and that $\theta_\B$ and $\theta_\A$ are the corresponding generators. This means that there must exist generators $\theta_1,\dots ,\theta_j\in D(\A)$ that combine to give $\theta_\B$ $$\theta_\B=\sum\limits_{i=1}^j p_i\theta_i .$$ But since $\B\subseteq \A$ we have that $\theta_i\in D(\B)$ for all $i=1,\dots, j$. This is a contradiction because all the $\theta_i$ have degree strictly less than $\theta_\B$ which would make it not a generator. \owari

\

Now we can state and prove the main theorem of this section.

\begin{theorem}\label{hypbound}

Suppose that $\A$ is a hypersolvable arrangement with hyperexponents $\mathrm{hypexp}(\A)=(b_1,\dots, b_k)$ and solvable filtration $\emptyset =\A_0\subset \A_1\subset \cdots \subset \A_k=\A$. Let $H\in \A\backslash \A_{k-1}$, put $\A'=\A\backslash H$, set $d_{\A'}$ and $d_\A$ to be the maximum degree of a set of generators of the module of derivations for $\A'$ and $\A$ respectively, and let $\rho(\A)=\max \{b_i \}_{i=1,\dots , k}$. Then $d_{\A'}\geq \rho(\A)-1$ and $d_\A\geq \rho(\A)-1$.
\end{theorem}

\noindent\textbf{Proof}. We induct on $k$. The base case $k=1$ is trivial because $d_{\A'}=0$ and $\rho(\A)-1=0$. Now assume $k>1$ and take $H\in \A\backslash \A_{k-1}$. By Corollary \ref{genbound} $d_{\A'}\geq |\A'|-|\A^H|=|\A'|-|\A_{k-1}|=|\A_k|-|\A_{k-1}|-1=b_k-1$ by Lemma \ref{restrictsolv}. Then by induction hypothesis we know that any deletion of the arrangements $\A_i$ for $1\leq i\leq k-1$ must have maximal degree generators bounded by $\rho(\A_j)$ for $1\leq j\leq k-1$. Hence by Lemma \ref{maxup} we are done. \owari

\

Unfortunately the lower bound found in Theorem \ref{hypbound} maybe far from tight for even simple examples. 

 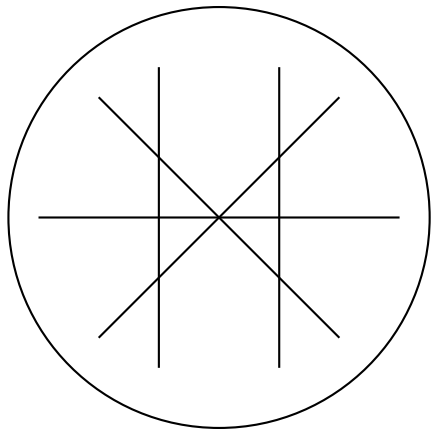
\begin{figure}\begin{center}\begin{tikzpicture}[scale=.8]
\draw[thick] (-2,-2)--(2,2);
\draw[thick] (2,-2)--(-2,2);
\draw[thick] (-1,-2.5)--(-1,2.5);
\draw[thick] (1,-2.5)--(1,2.5);
\draw[thick] (-3,0)--(3,0);

\draw[thick]  (0,0) circle (3.5cm);


\end{tikzpicture}\end{center}\caption{Hypersolvable not generic}\label{hypsol1ex}\end{figure}

\begin{example}\label{super1ex}

Let $\A$ be defined by the polynomial $Q(\A)=yz(x-z)(x+z)(x-y)(x+y)$ (see Figure \ref{hypsol1ex}). This example first appeared in \cite{Falk93} and was revisited in \cite{JP98}. A solvable filtration for $\A$ is $\emptyset\subset \{ x+z\}\subset \{x+z,z,x-z\} \subset \{x+z,z,x-z, x+y\}\subset \A$. Hence $\mathrm{hypexp}(\A)=\{1,2,1,2\}$. However a Macaulay2 (see \cite{M2}) calculation shows that the degrees of the generators of $D(\A)$ are $\{1,3,3,4\}$.
 So, $d_{\A}=4$ and $\rho(\A)-1=1$ which is far from tight for such a small example.
\end{example}

Clearly the study of derivations on hypersolvable arrangements is far from complete. Are there better bounds for $d_\A$? Can one use the rest of the hyperexponents to get bounds on the rest of the derivation degree sequence?







\bibliographystyle{amsplain}

\bibliography{non-freegens.bib}

\end{document}